\documentclass[12pt,psamsfonts]{amsart}

\newtheorem{anyprop}{Anyprop}[section]

\newtheorem{theorem}[anyprop]{Theorem}
\newtheorem{lemma}[anyprop]{Lemma}
\newtheorem{proposition}[anyprop]{Proposition}
\newtheorem{corollary}[anyprop]{Corollary}

\theoremstyle{definition}

\newtheorem{example}[anyprop]{Example}

\newtheorem{remark}[anyprop]{Remark}

\newtheorem{situation}[anyprop]{Situation}

\newtheorem{theoremintro}{Theorem}
\newtheorem{corollaryintro}{Corollary}

\newcommand{\NN}{\mathbb{N}}

\newcommand{\PP}{\mathbb{P}}

\newcommand  {\shM}     {\mathcal{M}}

\newcommand  {\shL}     {\mathcal{L}}

\newcommand  {\shS}     {\mathcal{S}}
\newcommand  {\shT}     {\mathcal{T}}

\newcommand  {\shQ}     {\mathcal{Q}}

\newcommand  {\foa}     {\mathfrak{a}}

\newcommand  {\fom}     {\mathfrak{m}}

\newcommand  {\dual}    {\vee}

\newcommand  {\lra}     {\longrightarrow}

\renewcommand{\O}       {\mathcal{O}}

\newcommand  {\Proj}    {\operatorname{Proj}}

\newcommand  {\ra}      {\rightarrow}

\newcommand  {\rk}    {\operatorname{rk}}

\newcommand  {\Syz}     {\operatorname{Syz}}

\newcommand{\komdots}{ , \ldots , }
\newcommand{\plusdots}{ + \ldots + }

\newcommand{\subsetdots}{ \subset \ldots \subset }

\newcommand{\numiii}{\renewcommand{\labelenumi}{(\roman{enumi})}}

\newcommand{\stacklra}[1]{ \stackrel{ #1 }{\lra} }

\newcommand{\length} {\lambda}

\theoremstyle{remark}

\numberwithin{equation}{section}

\usepackage{amscd}
\usepackage{amssymb}

\def\mydate{\number\day\space\ifcase\month \or January\or February\or March\or April\or May\or
June\or July\or August\or September\or October\or November\or
December\fi \space\number\year}

\begin{document}

\title[Hilbert-Kunz multiplicity]
{The rationality of the Hilbert-Kunz multiplicity in graded dimension two}

\author[Holger Brenner]{Holger Brenner}
\address{Department of Pure Mathematics, University of Sheffield,
  Hicks Building, Hounsfield Road, Sheffield S3 7RH, United Kingdom}
%\curraddr{}
\email{H.Brenner@sheffield.ac.uk}

%\date{}

% at present the "communicated by" line appears only in ERA and PROC
%\commby{}

%\dedicatory{Preliminary version,  \mydate}

\renewcommand{\arraystretch}{10}

\begin{abstract}
We show that the Hilbert-Kunz multiplicity is a rational number for
an $R_+-$primary homogeneous ideal $I=(f_1 \komdots f_n)$ in a
two-dimensional graded domain $R$ of finite type over an
algebraically closed field of positive characteristic. More
specific, we give a formula for the Hilbert-Kunz multiplicity in
terms of certain rational numbers coming from the strong
Harder-Narasimhan filtration of the syzygy bundle $\Syz(f_1 \komdots
f_n)$ on the projective curve $Y=\Proj R$.
\end{abstract}

\maketitle


%===========================================================
\section*{Introduction}

Suppose that $(R, \fom)$ is a local Noetherian Ring of dimension $d$
containing a field $K$ of positive characteristic $p$. Let $I=(f_1
\komdots f_n)$ denote an $\fom$-primary ideal, and denote by
$I^{[q]}=(f_1^{q} \komdots f_n^{q})$ the ideal given by the powers
$f^q$, where $q=p^{e}$. The ideal $I^{[q]}$ is the extended ideal of
$I$ under the $e$-th Frobenius homomorphism $R \ra R$, hence
independent of the choice of generators. Since $I$ is
$\fom$-primary, the length of the residue class ring $R/I^{[q]}$ is
finite for every prime power $q$.

The function $e \mapsto \length ( R/I^{[p^{e}]}) $, where $\length $ denotes the length,
is called the Hilbert-Kunz function
of the ideal $I$ and was introduced in \cite{kunzpositive}
(see also \cite{kunznoetherian}).
The Hilbert-Kunz multiplicity is defined as the limit
$$ e_{HK}(I)= \lim_{e \ra \infty}  \, \length (R/I^{[p^{e}]}) /p^{ed} \, .$$
This limit exists as a positive real number, as shown by Monsky in
\cite{monskyhilbertkunz} (see also \cite[Chapter
6]{hunekeapplication} and \cite[I.7.3]{robertsmultiplicity}). In the
same paper, Monsky writes, ``we suspect, but have no idea how to
prove, that $c(M)$[that is $e_{HK}$] is always rational''. C. Huneke
has put this question on his top ten list of problems in commutative
algebra \cite{huneketalksaltlakecity}. The Hilbert-Kunz multiplicity
of the maximal ideal is also called the Hilbert-Kunz multiplicity
$e_{HK}(R)$ of the ring itself and gives an important invariant, but
even in this case the rationality is only known in some special
cases.

Let us briefly recall what is known up to now. Most rationality
results deal only with the maximal ideal in hypersurface rings of
special type. Han and Monsky succeeded in the computation of the
Hilbert-Kunz multiplicity of a Brieskorn hypersurface
$X_0^{\delta_0} \plusdots X_N^{\delta_N}$ \cite{hanmonsky}, and
Conca provided a formula for it in the case of a homogeneous
binomial hypersurface of type $X_0^{\delta_0} \cdots X_m^{\delta_m}
- X_{m+1}^{\delta_{m+1}} \cdots X_N^{\delta_N}$
\cite{concahilbertkunz}. The rationality of the Hilbert-Kunz
multiplicity for cones $K[X,Y,Z]/(H)$ over a plane cubic curve
$V_+(H) \subset \PP^2$ was shown by Buchweitz, Chen and Pardue. It
is equal to $9/4$ in the smooth case and $7/3$ in the singular case,
see \cite{buchweitzchenhilbertkunz}. This was generalized by
Fakhruddin and Trivedi in \cite{fakhruddintrivedi} to any cone over
an elliptic curve. In his thesis \cite{teixeirathesis}, Teixeira
obtained the rationality of the Hilbert-Kunz multiplicity in the
case of hypersurfaces of type $H= \sum_i G_i(x_i,y_i)$, where the
$G_i$ are homogeneous.

For more general ideals, it is well known that the Hilbert-Kunz
multiplicity of a monomial ideal in a toric ring is rational
\cite{watanabetoric}. Watanabe and Yoshida give a formula for the
Hilbert-Kunz multiplicity of an integrally closed ideal in a
two-dimensional Gorenstein quotient singularity in terms of data
coming from a minimal resolution (see
\cite{watanabeyoshidahilbertkunzinequality},
\cite{watanabeyoshidahilbertkunztwo},
\cite{watanabeyoshidahilbertkunzmckay}). For other results
concerning estimates for the Hilbert-Kunz multiplicity and the
relationship to other ring-theoretic properties consult
\cite{blickleenescu}, \cite{haneshilbertkunz},
\cite{hunekeyaohilbertkunz},
\cite{millerfrobeniusprojectivedimension},
\cite{watanabeyoshidahilbertkunzminimal},
\cite{watanabeyoshidahilbertkunzthree}.

In this paper we show that the Hilbert-Kunz multiplicity $e_{HK}(I)$
is indeed a rational number, where $R$ is a normal standard-graded
two-dimensional domain of finite type over an algebraically closed
field of positive characteristic, and where $I$ denotes a
homogeneous $R_+$-primary ideal. The main idea is to use the short
exact sequence (set $d_i = \deg (f_i)$) of locally free sheaves
$$ 0 \lra \Syz (f_1^q \komdots f_n^q)(m) \lra \bigoplus_{i=1}^n \O(m-qd_i)
\stackrel{f_1^q \komdots f_n^q}{\lra} \O(m) \lra 0 \, $$
on the smooth projective curve $Y= \Proj R$
and to compute
$$ \length  ((R/I^{[q]})_m)
= h^0 (\O(m)) - \sum_{i=1}^n  h^0(\O(m-qd_i ))+ h^0(\Syz (f_1^q
\komdots f_n^q)(m) ) \, . $$ One computes $\length (R/I^{[q]})$ by
summing over $m$, which is a finite sum since this alternating sum
is $0$ for $m \gg 0$. Of course the crucial point is to control the
behavior of the global syzygies $H^0(Y, \Syz (f_1^q \komdots
f_n^q)(m))$ for different values of $q$ and $m$. Hence we are
concerned with a Frobenius-Riemann-Roch problem. We have $\Syz(f_1^q
\komdots f_n^q)(0) = F^{e*} ( \Syz(f_1 \komdots f_n )(0))$, where
$q=p^{e} $ and $F:Y \ra Y$ is the absolute Frobenius morphism; so we
have to understand the global sections of the Frobenius pull-backs
of the syzygy bundle $\Syz(f_1 \komdots f_n)$. The behavior of a
locally free sheaf under the Frobenius morphism is full of
surprising phenomena. For example, the Frobenius pull-back of a
semistable sheaf need not be semistable.

However, a recent theorem of A. Langer \cite[Theorem
2.7]{langersemistable} shows for a locally free sheaf $\shS$ on a
smooth projective variety $Y$ that the Harder-Narasimhan filtration
of some Frobenius pull-back has strongly semi\-stable quotients.
This means that for $e$ big enough there exists a filtration
$$ 0 = \shS_0 \subset \shS_1 \subsetdots \shS_t =F^{e*} (\shS) $$
such that the quotients
$ \shS_k/ \shS_{k-1}$ are strongly semistable
of decreasing slopes
$\mu_k( F^{e*}(\shS)) = \mu(\shS_k / \shS_{k-1})$.
Strongly semistable means that every
Frobenius pull-back is also semistable.

The existence of this strong Harder-Narasimhan filtration allows in
particular to define for $k=1 \komdots t$ rational numbers by
setting $\bar{\mu}_k (\shS) := \frac{\mu_k (F^{e}(\shS))}{ p^{e}} $.
These numbers, the length $t$ of the strong Harder-Narasimhan
filtration and the ranks $r_k= \rk (\shS_k/ \shS_{k-1})$ are all
independent of $q \gg 0$. Applying this to the syzygy bundle $\shS=
\Syz (f_1 \komdots f_n)$ we get numbers which control the global
syzygies for varying $q$. Therefore they enter (we set $ \nu_k = -
\bar{\mu}_k / \deg (Y)$) into the following simple formula for the
Hilbert-Kunz multiplicity, which is our main result and gives the
rationality (Theorem \ref{hilbertkunzformula}).

\begin{theoremintro}
Let $R$ denote a two-dimensional standard-graded normal domain
and let $I=(f_1 \komdots f_n)$ denote a homogeneous $R_+$-primary ideal
generated by homogeneous elements $f_i$ of degree $d_i, i=1 \komdots n$.
Then the Hilbert-Kunz multiplicity $e_{HK}(I)$ is given by
$$\frac{ \deg (Y)}{2}( \sum_{k=1}^t r_k \nu_k^2 - \sum_{i=1}^n d_i^2) \, .$$
In particular, the Hilbert-Kunz multiplicity is a rational number.
\end{theoremintro}

The rationality is also true without the assumptions normal and standard-graded
(see Corollary \ref{rationalcorollary}).
As an easy corollary we get a description for the Hilbert-Kunz multiplicity
for cones over plane curves (Corollary \ref{planecurve}).
This was only known for
degree $h \leq 3$ so far.
Furthermore, this corollary is independent of the result of Langer,
since the existence of the strong Harder-Narasimhan filtration
is clear by elementary means in this case.

\begin{corollaryintro}
Let $C=V_+(H) \subset \PP^2$ denote a smooth plane projective curve of degree $h$, $R=K[X,Y,Z]/(H)$.
Then there exists a rational number $ \frac{3}{2} \leq  \nu_2  \leq 2$
such that the Hilbert-Kunz multiplicity of $R$ is
$$e_{HK}(R) =h (\nu_2^2 -3 \nu_2 +3 ) \, .$$
\end{corollaryintro}

The rationality of the Hilbert-Kunz multiplicity for the maximal
ideal in dimension two was proved independently by V. Trivedi in
\cite{trivedihilbertkunz}. I thank M. Blickle and K. Watanabe for
useful discussions and the referee for useful comments.

\section{Preliminaries}
\label{pre}

We recall briefly some notions about vector bundles, see
\cite{huybrechtslehn} for details. Let $Y$ denote a smooth
projective curve over an algebraically closed field $K$. The degree
of a locally free sheaf $\shS$ on $Y$ of rank $r$ is defined by
$\deg (\shS) = \deg \bigwedge^r (\shS)$. The slope of $\shS$,
written $\mu(\shS)$, is defined by $\deg(\shS)/ r$. The slope has
the property that $\mu(\shS \otimes \shT)= \mu(\shS) + \mu(\shT)$.
Both the degree and the slope behave well under finite mappings: if
$\varphi: Y' \ra Y$ is a finite morphism of smooth projective curves
of degree $q$, then $\deg( \varphi^*(\shS)) =q \deg (\shS)$ and $\mu
(\varphi^* (\shS))= q \mu (\shS)$.

A locally free sheaf $\shS$ is called semistable if $\mu(\shT) \leq
\mu(\shS)$ holds for every locally free subsheaf $\shT \subseteq
\shS$. Dualizing and tensoring with an invertible sheaf does not
affect this property. For every locally free sheaf $\shS$ on $Y$
there exists the so-called Harder-Nara\-sim\-han filtration $\shS_1
\subsetdots \shS_t =\shS$, where the $\shS_k$ are locally free
subsheaves. This filtration is unique and has the property that the
quotients $\shS_{k}/ \shS_{k-1}$ are semistable and $\mu (\shS_{k}/
\shS_{k-1}) > \mu(\shS_{k+1}/ \shS_{k})$ holds.

The number $\mu(\shS_1) = \mu_{\rm max}(\shS)$ is called the maximal slope of $\shS$,
and the minimal slope of $\shS$ is
$ \mu_{\rm min}(\shS) =\mu(\shS/ \shS_{t-1})$.
The existence of global sections can be tested with the maximal slope: if $\mu_{\max}(\shS)<0$,
then $H^0(Y, \shS)=0$.
Furthermore we have the relation $\mu_{\max} (\shS)= - \mu_{\min}(\shS^\dual)$,
where $\shS^\dual$ denotes the dual bundle.

In positive characteristic, the pull-back of a semistable bundle
under the absolute Frobenius $F: Y \ra Y$ is in general not
semistable. If it stays semistable for every Frobenius power, then
the bundle is called strongly semistable, a notion introduced by
Miyaoka in \cite{miyaokachern}.

Consequently, the pull-back of the Harder-Narasimhan filtration of
$\shS$ under the Frobenius does not in general give the
Harder-Narasimhan filtration of $F^*(\shS)$. However, a recent
result of A. Langer \cite[Theorem 2.7]{langersemistable} shows that
there exists a Frobenius power $F^{e}$ such that the quotients in
the Harder-Narasimhan filtration of the pull-back $F^{e*}(\shS)$ are
all strongly semistable. We call such a filtration the strong
Harder-Narasimhan filtration and denote it by
$$ 0 \subset \shS_1^q \subsetdots \shS_t^q = F^{e*}(\shS) \, .$$
For $e' \geq e$ the Harder-Narasimhan filtration of
$F^{e' *} (\shS)$ is
$$ \shS_1^{q'}=F^{(e'-e)*}(\shS_1^q) \subsetdots \shS_t^{q'} = F^{(e'-e)*}(\shS_t^q) \, .$$
Using this we define rational numbers $ \bar{\mu}_k = \bar{\mu}_k
(\shS) = \frac{\mu(\shS_k^q/ \shS_{k-1}^q)}{q}$ for $q \gg 0$. The
length $t$ of the strong Harder-Narasimhan filtration as well as the
ranks $r_k= \rk( \shS^q_k/ \shS^q_{k-1})$ are independent of $q \gg
0$. For $q \gg 0$ we have $\mu (\shS_k^q/ \shS^q_{k-1}) =q
\bar{\mu}_k$ and $ \deg ( \shS^q_k/ \shS^q_{k-1} ) = q r_k
\bar{\mu}_k$. Furthermore note that
$$\rk (\shS_k^q) = r_1 \plusdots r_k \,\mbox{ and }\,
\deg (\shS_k^q) = q(r_1 \bar{\mu}_1 \plusdots r_k\bar{\mu}_k ) \,
.$$ We also set $\bar{\mu}_{\rm max} (\shS) = \bar{\mu}_1(\shS)$ and
$\bar{\mu}_{\rm min} (\shS) = \bar{\mu}_t(\shS)$

We shall apply these notions and facts to syzygy bundles. Let $R$
denote a normal standard-graded\footnote{Throughout this paper the
assumption standard-graded, meaning that $R$ is generated by
finitely many elements of degree one, might be weaken to the
property that $R$ is an $\NN$-graded domain of finite type and that
there exists finitely many elements $x_k$ of degree one such that
$\Proj R = \bigcup_k D_+(x_k)$. This last property is enough to
ensure that $\O(1)$ is an invertible sheaf and makes therefore
everything work (look e.g. at the proof of \cite[Proposition
II.5.12]{haralg}); see also Corollary \ref{rationalcorollary}
below.} domain over an algebraically closed field (of any
characteristic) and let $f_1 \komdots f_n$ denote homogeneous
generators of an $R_+$-primary ideal of degrees $d_1 \komdots d_n$.
This give rises to the short exact sequence of locally free sheaves
on $Y= \Proj R$,
$$ 0 \lra \Syz(f_1 \komdots f_n)(m) \lra
\bigoplus_{i=1}^n  \O(m-d_i) \stacklra {f_1 \komdots f_n}  \O(m)
\lra 0 \, .$$ For $m=0$ we write also $\Syz(f_1 \komdots f_n)$
instead of $\Syz(f_1 \komdots f_n)(0)$. Due to this defining
sequence, the rank of the syzygy bundle is $n-1$ and its degree is
$((n-1)m- \sum_{i=1}^n d_i) \deg \O_Y(1)$.

Now suppose that the algebraically closed ground field $K$
has positive characteristic.
The pull-back of the short exact sequence of locally free sheaves
under the $e$-th absolute Frobenius morphism $F^{e}: Y \ra Y$ yields
$$ 0 \lra (F^{e}(\Syz(f_1 \komdots f_n))) (m) \lra
\bigoplus_{i=1}^n  \O(m-qd_i) \stacklra {f_1^q \komdots f_n^q}  \O(m) \lra 0 \, $$
(pull-back the sequence for $m=0$ and tensor it with $\O(m)$ again).
Therefore $(F^{e}(\Syz(f_1 \komdots f_n))) (m) = \Syz(f_1^q \komdots f_n^q)(m)$.
We want to compute $\length ((R/I^{[q]} )_m)$ using
this exact sequence.
The global sections $\Gamma(Y, -)$ of this sequence
yield (since $R$ is assumed to be normal and standard-graded)
$$ 0 \lra \Gamma(Y, \Syz (f_1^q \komdots f_n^q)(m))
\lra \bigoplus_{i=1}^n R_{ m-qd_i} \stackrel{f_1^q \komdots f_n^q}
{\lra} R_m  \lra \ldots $$ and the cokernel of the last mapping is
$(R /I^{[q]}) _m$. This is the same as the kernel of the mapping
$H^1( Y, \Syz(f_1^q \komdots f_n^q)(m)) \ra \bigoplus_{i=1}^n H^1(Y,
\O(m-qd_i))$. Hence we compute
$$ \length  ((R/I^{[q]}) _m)
= h^0 (\O(m)) - \sum_{i=1}^n  h^0(\O(m-qd_i ))+ h^0(\Syz (f_1^q \komdots f_n^q)(m) )  \, $$
and then we sum over $m$.

%If the syzygy bundle $\Syz(f_1^q \komdots f_n^q)(m)$ is strongly semistable,
%then its maximal slope is negative
%for $m  <  q \sum d_i /(n-1)$. In this case
%there exist no global non-trivial syzygy and hence we have the injection
%$0 \lra \bigoplus_{i=1}^n H^0(Y, \O( m-qd_i))  \lra H^0(Y, \O(m)) $
%and therefore in this range the global syzygies do not enter the computation.

\section{The case of a strongly semistable syzygy bundle}

In this section we prove some results
about the asymptotic behavior of $h^0(\O(m-qd_i ))$ and $h^0(\O(m))$ and we
apply this to compute the Hilbert-Kunz multiplicity
under the condition that the syzygy bundle is strongly semistable.
We fix the following situation.

\begin{situation}
\label{situation1} Let $R$ denote a normal two-dimensional
standard-gra\-ded domain over an algebraically closed field of
positive characteristic $p$ with corresponding smooth projective
curve $Y= \Proj R$ of genus $g$. Set $\deg (Y) = \deg (\O_Y(1))$.
Let $(f_1 \komdots f_n)$ denote a homogeneous, $R_+$-primary ideal
given by homogeneous ideal generators of degree $d_i$. Let $q=p^{e}$
denote varying prime powers.
\end{situation}

We will often use the notation $O(g)$ for the asymptotic behavior of
a function $f(q)$ in one variable. The equation $f= O(g)$ means that
$f/g$ is bounded for $q \ra \infty$. The functions we consider will
be defined only for prime powers $q=p^{e}$, hence $f=O(q)$ means
that $f(q)/q$ is bounded. Such functions are negligible in our
situation, since then $f(q)/q^2 \ra 0$.

\begin{lemma}
\label{invertible}
Suppose the situation of \ref{situation1}.
Let $\nu$ denote a positive rational number.
For $i=1 \komdots n$, if $\nu \geq d_i$, we have
$$ \sum_{m=0}^{\lceil q \nu \rceil} h^0(\O(m-qd_i))
= q^2 \frac{\deg (Y)}{2} (  \nu-d_i)^2      + O(q) \, .$$
and
$$ \sum_{m=0}^{\lceil q \nu \rceil} h^0(\O(m))
= q^2 \frac{\deg (Y)}{2}  \nu ^2       + O(q) \, .$$
\end{lemma}
\proof
By Riemann-Roch we have
$$h^0(\O(m-qd_i)) =( m-qd_i) \deg (Y) +1-g +h^1(\O(m-qd_i)) \, .$$
Therefore we have $\sum_{m=0}^{\lceil q \nu \rceil}
h^0(\O(m-qd_i))=$
\begin{eqnarray*}
&=& \sum_{m=qd_i}^{\lceil q \nu \rceil} h^0(\O(m-qd_i)) \cr &=&
\sum_{m=qd_i}^{\lceil q \nu \rceil} \big( ( m-qd_i) \deg (Y) +1-g
+h^1(\O(m-qd_i))   \big) \cr &=&  \deg (Y) \! \!
\sum_{m=qd_i}^{\lceil q \nu \rceil} ( m-qd_i) + (\lceil q \nu \rceil
- qd_i+1 )(1-g) + \! \!   \sum_{m=qd_i}^{\lceil q \nu \rceil} \! \!
h^1(\O(m-qd_i)) \cr &=& \frac{\deg (Y)}{2}(\lceil q \nu \rceil -
qd_i +1 )(\lceil q \nu \rceil - qd_i) + O(q) +O(q^0) \cr &=&
\frac{\deg (Y)}{2}q^2 ( \nu - d_i )^2 + O(q)
\end{eqnarray*}
Here we used on the right that $H^1(Y, \O(m-qd_i))=0$ for $m-qd_i \gg 0$, and this bound is independent of $q$.
The proof for the statement about $\O(m)$ is the same.
\qed

\begin{lemma}
\label{minslopelemma}
Let $\shS$ denote a locally free sheaf on a smooth projective curve $Y$ with a very ample invertible sheaf $\O(1)$ of degree $\deg( \O(1))=\deg(Y)$. Denote the pull-back of $\shS$ under the $e$-th
absolute Frobenius by
$\shS^q$, $q=p^{e}$.
Set $\nu= - \frac{ \bar{\mu}_{\min} (\shS)} {\deg (Y)}$.
Then
$$\sum_{m= \lceil q \nu \rceil} ^ \infty
h^1(Y,  \shS^q(m)) = O(q) \, .$$
\end{lemma}
\proof The minimal slope of $\shS^q$ is $- q \deg (Y) \nu$ for $q$
big enough. By Serre duality we have
$$h^1(\shS^q(m))
= h^0((\shS^{q})^ \dual (-m) \otimes \omega_Y) \, .$$ Now for $m
> \lceil q \nu  \rceil  +\deg( \omega_Y)/\deg (Y)$ we have
\begin{eqnarray*}
\mu_{\max} (( \shS^q )^\dual(-m) \otimes \omega_Y)
&=&- \mu_{\min}( \shS^q (m)) + \mu( \omega_Y)  \cr
&=& - \mu_{\min}( \shS^q ) - m \deg (Y) + \deg( \omega_Y) \cr
&=& -q \bar{\mu}_{\min} (\shS) - m \deg (Y) + \deg (\omega_Y) \cr
& <& 0 \, .
\end{eqnarray*}
So for these $m$ we have $H^1(Y,  \shS^q (m) )=0$
and our sum is indeed finite
running in the range
$ \lceil q \nu \rceil \leq m \leq \lceil q \nu \rceil + \frac{\deg (\omega_Y)}{ \deg(Y)}$.
In particular, the length of this range is independent of $q$.

There exists a surjection $\bigoplus_{j \in J} \O( \alpha_j) \ra
\shS \ra 0$. Pulling this back under $F^e$ we get surjections
$\bigoplus_{j \in J} \O( q\alpha_j) \ra \shS^q \ra 0$ and therefore
$$ \bigoplus_{j \in J} H^1(Y, \O(q \alpha_j + m)) \lra H^1(Y,\shS^q (m)) \lra 0\, .$$

For $m$ fulfilling $\lceil q \nu  \rceil \leq m \leq \lceil q \nu
\rceil  + \frac{\deg( \omega_Y)}{\deg (Y)}$ we see that $q \alpha_j
+ m$ varies in a range between $ q \beta_j $ and $q \beta_j +c$ ($c$
and $\beta_j$ independent of $q$). We have to understand the
asymptotic behavior of $h^1(\O(q \beta_j + \ell))$, $0 \leq \ell
\leq c$ for $q$ large. But $h^1(\O(q \beta_j + \ell)) = h^0(\O(-q
\beta_j - \ell ) \otimes \omega_Y)$ goes to $0$ for $ \beta_j
>0$ and it is $O(q)$ for $ \beta_j \leq 0$. So in any case the first
cohomology is $O(q)$ and the same is true for the finite sums over
all $\ell$ and $j \in J$. \qed

\begin{corollary}
\label{strongsemistablelemma}
Suppose the situation \ref{situation1}.
Suppose that the syzygy bundle $\Syz(f_1 \komdots f_n)$ is strongly semi\-stable.
Then
$$\sum_{m= \lceil q \frac{d_1 \plusdots d_n}{n-1} \rceil}^\infty
h^1(\Syz(f_1^q \komdots f_n^q)(m)) = O(q) \, .$$
\end{corollary}
\proof
This follows from Lemma \ref{minslopelemma} applied to $\shS=\Syz(f_1 \komdots f_n)$,
since in this case $\nu = \frac{d_1 \plusdots d_n}{n-1}$.
\qed

\begin{theorem}
\label{strongsemistabletheorem}
Suppose the situation of \ref{situation1}.
Suppose that the syzygy bundle is strongly semistable.
Then the Hilbert-Kunz function of $I=(f_1 \komdots f_n)$ may be written as
$$ \varphi(I,q) =  q^2 \frac{\deg (Y)}{2}
(\frac{(\sum_i d_i)^2}{n-1}  - \sum_i d_i^2)
+O(q)$$
and the Hilbert-Kunz multiplicity is
$$e_{HK}(I)=\frac{\deg (Y)}{2}
( \frac{(\sum_i d_i)^2}{n-1}  - \sum_i d_i^2) \, .$$
In particular it is a rational number.
\end{theorem}
\proof
We have the equations $\length (R/I^{[q]}) =$
\begin{eqnarray*}
&=& \sum_{m=0}^\infty \length ((R/ I^{[q]})_m) \cr
&=& \sum_{m=0}^{\lceil q \frac{d_1 \plusdots d_n}{n-1} \rceil -1}  \length ((R/ I^{[q]}) _m) + O(q)
 \cr
&=& \sum_{m=0}^{\lceil q \frac{d_1 \plusdots d_n}{n-1} \rceil -1} h^0(\O(m))
-\sum_{i=1}^n
\big(\sum_{m=0}^{\lceil q \frac{d_1 \plusdots d_n}{n-1} \rceil -1} h^0(\O(m-qd_i))\big) + O(q) \cr
&=&q^2 \frac{\deg (Y)}{2}
\big((\frac{d_1 \plusdots d_n}{n-1})^2 - \sum_{i=1}^n ( \frac{d_1 \plusdots d_n}{n-1} - d_i )^2     \big)
+O(q) \, .
\end{eqnarray*}
Here the second equation is due to Corollary
\ref{strongsemistablelemma}, since we have the inclusion
$(R/I^{[q]})_m \subseteq H^1(Y, \Syz(f_1^q \komdots f_n^q)(m))$ and
so the sum over $m \geq \lceil q \frac{d_1 \plusdots d_n}{n-1}
\rceil$ behaves like $O(q)$. The third equation is due to the fact
that $\Syz (f_1^q \komdots f_n^q)(m)$ is semistable and of negative
degree in the given range $m \leq \lceil  q \frac{d_1 \plusdots
d_n}{n-1}\rceil -1$, hence it has no global non-trivial sections.
The fourth equation is due to Lemma \ref{invertible} applied for
$\nu = \frac{d_1 \plusdots d_n}{n-1}$ (we could have written $\lceil
q \nu \rceil -1$ instead of $\lceil q \nu \rceil$ in the lemma).

Up to the factor $\deg (Y)/2$ we may write the Hilbert-Kunz
multiplicity as
\begin{eqnarray*}
(\frac{\sum_j d_j}{n-1})^2 - \sum_{i=1}^n (\frac{\sum_j d_j}{n-1}-d_i)^2
&=&  (\frac{\sum_j d_j}{n-1})^2
- \sum_{i=1}^n \big( (\frac{\sum_j d_j}{n-1})^2 +d_i^2 -2d_i \frac{\sum_j d_j}{n-1} \big) \cr
&=& -(n-1) (\frac{\sum_j d_j}{n-1})^2   - \sum_i d_i^2 +2 \frac{(\sum_j d_j)^2}{n-1} \cr
&=& \frac{(\sum_i d_i)^2}{n-1}  - \sum_i d_i^2  \,
\end{eqnarray*}
\qed

\begin{corollary}
\label{strongsemistableconstant}
Suppose the situation of \ref{situation1}.
Suppose that the syzygy bun\-dle is strongly semistable
and that the degrees of the ideal generators are con\-stant, $d=d_i$.
Then the Hilbert-Kunz multiplicity of $I=(f_1 \komdots f_n)$ is
$$e_{HK}(I)=\frac{\deg (Y)}{2}\frac{nd^2}{n-1} \, .$$
\end{corollary}
\proof
This follows from Theorem \ref{strongsemistabletheorem}, since
$$\frac{(nd)^2}{n-1} -nd^2 = \frac{n^2d^2-n(n-1)d^2}{n-1} = \frac{n}{n-1} d^2\, . $$\qed

\begin{corollary}
Let $Y \subset \PP^N$ denote a smooth projective curve, $Y= \Proj
R$, $R=K[X_0 \komdots X_N]/ \foa$. Suppose that the restriction of
the tangent bundle $\shT_{\PP^N}$ to the curve is strongly
semistable. Then the Hilbert-Kunz multiplicity of $R$ is $e_{HK}(R)=
\frac{\deg(Y)}{2} \frac{N+1}{N}$.
\end{corollary}
\proof We have to compute the Hilbert-Kunz multiplicity of the
maximal ideal $\fom=(X_0 \komdots X_N)$. The syzygy bundle $\Syz(X_0
\komdots X_N)$ on $\PP^N$ is the same as the cotangent bundle of
$\PP^N$ due to the Euler sequence (see \cite[Theorem II
8.13]{haralg}). Hence the strong semistability of the restriction of
this bundle implies the statement by Corollary
\ref{strongsemistableconstant} with $d=1$, $n=N+1$. \qed

\begin{example}
Consider the case of the maximal ideal $\fom =(X,Y,Z)$ on a smooth plane curve
$C= \Proj R$, $R= K[X,Y,Z]/(H)$ of degree $h= \deg (H) = \deg (C)$.
If the syzygy bundle $\Syz(X,Y,Z)|C$ is strongly semi\-sta\-ble,
then we get $e_{HK}(R)= \frac{3h}{4}$.
It was indeed a result of \cite[Corollary 1]{buchweitzchenhilbertkunz} that there
exists for every degree $h \geq 2$ a plane curve of degree $h$ whose Hilbert-Kunz multiplicity
is $\frac{3h}{4}$.
\end{example}

\section{Main results}

We treat now the case of an arbitrary syzygy bundle on a curve
making use of the strong Harder-Narasimhan filtration. The knowledge
of the Harder-Narasimhan filtration of a locally free sheaf $\shS$
contains a lot of information about the behavior of the global
sections $H^0(Y, \shS(m))$, as the following Lemma shows.

\begin{lemma}
\label{globalsection}
Let $\shS$ denote a locally free sheaf on a smooth projective curve $Y$ of genus $g$
over an algebraically closed field.
Let $\O_Y(1)$ be a very ample invertible sheaf and set
$\deg (Y) = \deg \O_Y(1)$.
Let $\shS_1 \subsetdots \shS_t =\shS$ be the Harder-Narasimhan filtration of $\shS$
and let $\mu_k(\shS) = \mu (\shS_k/ \shS_{k-1})$ denote the slopes of the semistable quotient
sheaves in this filtration and set $r_k= \rk (\shS_k/\shS_{k-1})$.
Then we have the following statements about the global sections and the first cohomology
of $\shS(m)$.

\numiii

\begin{enumerate}

\item
For $\mu_1 (\shS(m)) <0$ we have $H^0(Y, \shS(m))=0 $

\item
Fix $k$, $1 \leq k \leq t-1$. Let $m$ be such that $\mu_k (\shS (m))
> \deg (\omega_Y) $ and $ \mu_{k+1} (\shS(m)) <0 $.
Then  $H^0(Y,\shS(m)) \cong H^0(Y, \shS_k(m))$ and
\begin{eqnarray*}
h^0(\shS(m)) &=& \deg ( \shS_k(m))   + \rk (\shS_k) (1-g) \cr & =&
m(r_1 \plusdots r_k) \deg (Y) + r_1 \mu_1 \plusdots r_k \mu_k+ \rk
(\shS_k) (1-g)
\end{eqnarray*}

\item
For $\mu_t( \shS(m)) > \deg (\omega_Y) $ we have $H^1(Y,\shS(m) ) =0$.
\end{enumerate}
\end{lemma}
\proof The condition in (i) means that the maximal slope of the
sheaf $\shS(m)$ is negative, hence it cannot have non-trivial global
sections. (iii). By Serre duality we have $h^1(\shS(m) ) =
h^0(\shS^{\dual} (-m) \otimes \omega_Y)$ and we have
$$\mu_{\max} (\shS^{\dual} (-m) \otimes \omega_Y) =- \mu_{\min} ( \shS(m)) + \deg (\omega_Y)
= - \mu_t( \shS(m)) + \deg (\omega_Y) <0 \, ,
$$
hence $H^1(Y, \shS(m)) = 0$.

(ii).
Look at the sequence
$0 \ra \shS_k(m) \ra \shS(m) \ra \shS(m)/\shS_{k}(m) =\shQ_k(m)$.
Then we have $ \mu_{\min} (\shS_k (m)) = \mu ( (\shS_k/\shS_{k-1} )(m)) = \mu_k(\shS(m))> \deg (\omega_Y) $
and $ \mu_{\max} (\shQ_k (m)) = \mu( (\shS_{k+1}/ \shS_k)(m)) = \mu_{k+1} (\shS(m)) < 0$.
Due to this last observation,
$\shQ_k(m)$ does not have global non-trivial sections.
Therefore
$$H^0(Y, \shS(m)) \cong H^0(Y, \shS_k(m)) \, .$$
By Riemann Roch for locally free sheaves we have
$$h^0(\shS_k(m)) = \deg ( \shS_k(m))  + \rk (\shS_k)(1-g)  +h^1(\shS_k(m)) \, .$$
Serre duality yields $h^1(\shS_k(m)) = h^0((\shS_k(m))^\dual \otimes
\omega_Y)$. Now
$$
\mu_{\max} ((\shS_k(m))^\dual \otimes \omega_Y)= -\mu_{min} (\shS_k (m)) + \deg (\omega_Y) < 0 \, .$$
Hence $(\shS_k(m))^\dual \otimes \omega_Y$ does not have any non-trivial global section,
therefore $h^1(\shS_k(m))=0$
and we obtain the first equation.
The second equation is clear due
to
\begin{eqnarray*}
\deg (S_k(m)) &=& \deg (S_k) + \rk (\shS_k) \deg (\O(m)) \cr
&=&  r_1\mu_1 \plusdots r_k \mu_k + m(r_1 \plusdots r_k) \deg (Y) \, .
\end{eqnarray*}
\qed

\medskip
We fix now the situation and notation for the results in this section.

\begin{situation}
\label{situation} Let $R$ denote a two-dimensional normal
standard-graded (see footnote 1) domain over an algebraically closed
field $K$ of positive characteristic $p$ and let $Y= \Proj R$ denote
the corresponding smooth projective curve of genus $g$. Let $\deg
(Y)$ denote the degree of $\O_Y(1)=\O(1)$ given by $R$. Let $I=(f_1
\komdots f_n)$ denote an $R_+$-primary homogeneous ideal generated
by homogeneous elements $f_i$ of degree $d_i$. Let $\Syz(f_1
\komdots f_n)$ be the syzygy bundle on $Y$ and let $\bar{\mu}_k$
{\rm(}$r_k${\rm)} denote the slopes {\rm(}ranks{\rm)} of the
quotients in the strong Harder-Narasimhan filtration as explained in
the preliminaries and let $t$ denote its length. It is convenient to
set $\nu_k := - \bar{\mu}_k/ \deg (Y)$.
\end{situation}

\begin{remark}
The sum $r_1 \plusdots r_t= n-1$ equals the rank of the syzygy bundle $\Syz(f_1 \komdots f_n)$.
We have the relationship
$\sum_{k=1}^t \nu_k r_k = \sum_{i=1}^n d_i$
and the estimates
$$ \min_i (d_i)  \leq \nu_1 < \ldots < \nu_t \leq  \max_{i \neq j} (d_i +d_j) \, .$$
This follows from the defining sequence for the syzygy bundle and its Koszul resolution.
If $t=1$, then $\nu_1= \frac{d_1 \plusdots d_n}{n-1}$.
Think of the rational numbers $\nu_k$ as degree thresholds, where
something happens in the behavior of the global syzygies
$H^0(Y,\Syz(f_1^q \komdots f_n^q)(m))$,
when $m$ passes through $\nu_k$.
\end{remark}

\begin{proposition}
\label{globalsyzygy}
Suppose the situation and the notation of \ref{situation}.
Let $e$ be big enough such that
the Harder-Narasimhan filtration
of $F^{e}(\Syz(f_1 \komdots f_n)) =\Syz(f_1^q \komdots f_n^q)$ is strong, $q=p^{e}$.
Then the global syzygies have the following description.

\numiii

\begin{enumerate}

\item
For $m < q\nu_1 =-q \frac{ \bar{\mu}_1}{\deg (Y)}
=-q \frac{\bar{\mu}_{\max}(\Syz(f_1 \komdots f_n))}{\deg (Y)} $
we have
$$ H^0(Y, \Syz(f_1^q \komdots f_n^q) (m))  =0 \, .$$

\item
Fix $k$, $1 \leq k \leq t-1$.
For
$ q \nu_k + \frac{\deg (\omega_Y)}{\deg (Y)}
 < m  < q \nu_{k+1} $
we have
$$ h^0(Y, \Syz(f_i^q )  (m))
\! = \! \! q(r_1 \bar{\mu}_1 \plusdots r_k \bar{\mu}_k )+ m (r_1 \plusdots r_k) \deg (Y) + \rk (\shS_k)(1-g)\, .$$

\item
For
$m > q \nu_t + \frac{\deg (\omega_Y)}{\deg (Y)}
 =-q \frac{\bar{\mu}_t}{\deg (Y)} + \frac{\deg (\omega_Y)}{\deg (Y)}
$
we have $$H^1(Y, \Syz(f_1^q \komdots f_n^q)(m) )=0 \, .$$
\end{enumerate}
\end{proposition}

\proof
This follows from Lemma \ref{globalsection} applied to $\shS= \Syz(f_1^q \komdots f_n^q)$.
One only has to observe that
\begin{eqnarray*}
\bar{\mu}_{k} (\Syz(f_1^q \komdots f_n^q)(m))
&=& \bar{\mu}_{k} (\Syz(f_1^q \komdots f_n^q)) + \deg (\O(m)) \cr
&=& q \bar{\mu}_{k} + m \deg (Y)  \, .
\end{eqnarray*}
So for example the condition in the Lemma, that
$ \mu_k (\Syz(f_1^q \komdots f_n^q)(m)) > \deg (\omega_Y)$,
is equivalent with
$m \deg(Y) + q \bar{\mu}_k > \deg (\omega_Y)$
and hence with
$m > -q  \frac{\bar{\mu}_k}{\deg (Y)} + \frac{\deg(\omega_Y)}{\deg(Y)}$.
\qed

\begin{proposition}
\label{inbetween}
Suppose the situation and notation of \ref{situation}.
Fix $k$, $1 \leq k \leq t-1$.
Then we have
$ \sum_{m= \lceil q \nu_k \rceil} ^{ \lceil q \nu_{k+1} \rceil -1 } h^0(\Syz(f_1^q \komdots f_n^q)(m)) =$\,
$$q^2\deg (Y) \big( \frac{\nu_{k+1}^2- \nu_k^2}{2} (r_1 \plusdots r_k)
- (\nu_{k+1} -\nu_{k})(r_1 \nu_1 \plusdots r_k \nu_k)\big)
+ O(q) \, .
$$
\end{proposition}
\proof
We may assume that $q$ is big enough such that the
Harder-Narasimhan filtration of $\Syz(f_1^q \komdots f_n^q)$ is strong.
Let $\shS_k^q$ denote the $k$-th subbundle in the Harder-Narasimhan filtration.
We have
\begin{eqnarray*}
& & \sum_{m= \lceil  q \nu_k \rceil} ^{ \lceil  q \nu_{k+1} \rceil -1}
h^0(\Syz(f_1^q \komdots f_n^q)(m)) \cr
&=& \sum_{m= \lceil  q \nu_k \rceil} ^{ \lceil  q \nu_{k+1} \rceil -1}
h^0(\shS_k^q(m)) \cr
\vspace{1000pt}
&=&  \sum_{m= \lceil  q \nu_k \rceil} ^{\lceil  q \nu_{k+1} \rceil -1}
\big( \deg( \shS_k^q) + m \rk (\shS_k^q) \deg (Y)
+\rk (\shS_k^q) (1-g) +h^1(\shS_k^q (m))  \big) \cr
&=&  \sum_{m= \lceil  q \nu_k \rceil} ^{\lceil  q \nu_{k+1} \rceil -1}
\big( \deg( \shS_k^q) + m \rk (\shS_k^q) \deg (Y) \big)
+O(q)
+ \sum_{m=\lceil  q \nu_k \rceil}
^{\lceil  q \nu_k \rceil +\frac{ \deg (\omega_Y)}{ \deg(Y)} }
h^1(\shS_k^q (m)) \cr
&=&  \sum_{m= \lceil  q \nu_k \rceil} ^{\lceil  q \nu_{k+1} \rceil -1}
\big(q (r_1 \bar{\mu}_1 \plusdots r_k \bar{\mu}_k) + m \rk (\shS_k^q) \deg (Y) \big)
+ O(q) \cr
&=& \deg (Y) \big(\! \sum_{m= \lceil   q\nu_k \rceil} ^{\lceil  q \nu_{k+1}  \rceil -1}
 m (r_1 \plusdots r_k)
-q \! \! \! \! \sum_{m= \lceil   q\nu_k \rceil} ^{\lceil  q
\nu_{k+1} \rceil -1 } ( r_1 \nu_1 \plusdots r_k \nu_k ) \big)
+O(q) \cr &= & \deg (Y) \big(\frac{\lceil q\nu_{k+1} \rceil ( \lceil
q \nu_{k+1} \rceil - 1)}{2} - \frac{( \lceil q \nu_k \rceil +1)
\lceil q \nu_k \rceil)}{2} \big) (r_1 \plusdots r_k)\cr & &
\hspace{4cm} -q( \lceil q\nu_{k+1}\rceil - \lceil q\nu_k \rceil)
(r_1 \nu_1 \plusdots r_k \nu_k) + O(q)) \cr &=& \deg (Y) \big(q^2
\frac{ \nu_{k+1}^2 - \nu_k^2}{2} (r_1 \plusdots r_k)
%\cr & & \hspace{5cm}
-q^2 (\nu_{k+1} - \nu_k)(r_1 \nu_1 \plusdots r_k
\nu_k)\big) + O(q)
\end{eqnarray*}
The first equation is due to Lemma \ref{globalsection} (ii) applied
for $\shS = \Syz(f_1^q \komdots f_n^q)$. The second equation is due
to Riemann-Roch. The third equation was established in Proposition
\ref{globalsyzygy}(ii) (see also the proof of Lemma
\ref{globalsection}(ii)). For the fourth equation look at Lemma
\ref{minslopelemma}. Fix $\bar{q}$ big enough such that the
Harder-Narasimhan filtration is strong. Then $\bar{\mu}_{\min}
(\shS_k^{\bar{q}}) = \mu ( \shS_k^{ \bar{q}} / \shS_{k-1}^{\bar{q}})
= \bar{q} \bar{\mu}_k $. Applying Lemma \ref{minslopelemma} to $\shS
= \shS_k^{\bar{q}}$ with $\nu= \bar{q} \nu_k$ we get
$$\sum_{m= \lceil \tilde{q} (\bar{q} \nu_k) \rceil}^\infty
h^1((\shS^{\bar{q}}_k)^{\tilde{q}}(m)) = O( \tilde{q})  \, .$$ So
for $q= \tilde{q}\bar{q}$ we get the needed result, since $\tilde{q}
\leq q$. The other equations are clear. \qed

\medskip
We come now to the main result of this paper.

\begin{theorem}
\label{hilbertkunzformula}
Suppose the situation \ref{situation}.
The Hilbert-Kunz multiplicity $e_{HK}(I)$ is given by
$$\frac{ \deg (Y)}{2}( \sum_{k=1}^t r_k \nu_k^2 - \sum_{i=1}^n d_i^2) \, .$$
In particular, the Hilbert-Kunz multiplicity is a rational number.
\end{theorem}

\proof
We compute the length as $\length ( R/I^{[q]})=$
\renewcommand{\arraystretch}{10}
\begin{eqnarray*}
&=& \sum_{m=0}^{ \infty } \length (( R/I^{[q]})_m) \cr &=&
\sum_{m=0}^{ \lceil q \nu_t  \rceil -1 } \length (( R/I^{[q]})_m)  +
O(q) \cr &=&  \sum_{m=0}^{ \lceil q \nu_t \rceil -1} \big( h^0 (
\O(m))- \sum_{i=1}^n h^0(\O(m-qd_i))+h^0(\Syz(f_j^q )(m))\big)+O(q)
\cr &=& \sum_{m=0}^{  \lceil q \nu_t  \rceil -1 } h^0 ( \O(m)) -
\sum_{i=1}^n \big(   \sum_{m=0}^{  \lceil q \nu_t \rceil -1}
h^0(\O(m-qd_i)) \big) \cr & & \hspace{6.2cm} +\sum_{m=0}^{ \lceil q
\nu_t  \rceil -1} h^0(\Syz(f_j^q )(m)) + O(q) \cr &=& q^2
\frac{\deg(Y)}{2} \big( \nu_t^2- \sum_{i=1}^n (\nu_t -d_i)^2 \big)
+\sum_{m=0}^{ \lceil q \nu_t \rceil -1}h^0(\Syz(f_j^q )(m))+O(q)\, .
\end{eqnarray*}
Here the second equation is due to Lemma \ref{minslopelemma} (as in
the proof of Theorem \ref{strongsemistabletheorem}) and the last
equation is due to Lemma \ref{invertible}. We may write the term on
the right $\sum_{m=0}^{\lceil q \nu_t \rceil -1}h^0(\Syz(f_j^q)(m))
=$
\begin{eqnarray*}
&=& \sum_{m=0}^ {\lceil q \nu_1 \rceil -1} h^0(\Syz(f_j^q )(m))
+\sum_{k=1}^{t-1}
 \big( \sum_{m= \lceil q \nu_k \rceil} ^{\lceil q\nu_{k+1} \rceil -1}h^0(\Syz(f_j^q )(m))\big)
\cr
&=& \sum_{k=1}^{t-1}
 \big( \sum_{m= \lceil q \nu_k \rceil} ^{\lceil q\nu_{k+1} \rceil -1}  h^0(\Syz(f_j^q )(m))\big)
\cr
&=& q^2 \deg (Y) \sum_{k=1}^{t-1}
\big(\frac{ \nu_{k+1}^2 - \nu_k^2}{2} (r_1 \plusdots r_k) \cr
& & \hspace{5cm} -(\nu_{k+1} - \nu_k)(r_1 \nu_1 \plusdots r_k \nu_k)\big) + O(q)    \, .
\end{eqnarray*}
Here the second equation is due to
Proposition \ref{globalsyzygy} (i) and the last equation is due to Proposition \ref{inbetween}.
This yields for the Hilbert-Kunz multiplicity $e_{HK}(I)$ the expression
$$\frac{ \deg (Y)}{2}
\big( \nu_t^2 - \sum_{i=1}^n (\nu_t -d_i)^2
+  \sum_{k=1}^{t-1}       (\nu_{k+1}^2 - \nu_k^2)(r_1 \plusdots r_k) \hspace{2cm} $$
$$ \hspace{5cm} - 2 \sum_{k=1}^{t-1}
( \nu_{k+1} - \nu_{k})(r_1 \nu_1 \plusdots r_k \nu_k) \big) \, .$$
Now we can simplify. We have
$$ \nu_t^2 - \sum_{i=1}^n (\nu_t -d_i)^2
= - (n-1) \nu_t^2+ 2 \nu_t \sum_{i=1}^n d_i - \sum_{i=1}^n d_i^2 \, .$$
Furthermore we have
$$ \sum_{k=1}^{t-1} (\nu_{k+1}^2-\nu_k^2)(r_1 \plusdots r_k)
= - \sum_{k=1}^{t-1} r_k \nu_k^2 + \nu_t^2(r_1 \plusdots r_{t-1})$$
and similarly
$$ \sum_{k=1}^{t-1} (\nu_{k+1}-\nu_k)(r_1\nu_1 \plusdots r_k\nu_k)
= - \sum_{k=1}^{t-1} r_k\nu_k^2  + \nu_t(r_1\nu_1 \plusdots r_{t-1}\nu_{t-1}) \, .$$
Using the relations $r_1 \plusdots r_{t-1} =n-1-r_t$
and $r_1\nu_1 \plusdots r_{t-1}\nu_{t-1} = \sum_{i=1}^n d_i - r_t \nu_t$
we get alltogether
\begin{eqnarray*}
& & - (n-1) \nu_t^2+ 2 \nu_t \sum_{i=1}^n d_i - \sum_{i=1}^n d_i^2
-  \sum_{k=1}^{t-1} r_k\nu_k^2 + \nu_t^2(n-1-r_t) \cr
& & -2\big(- \sum_{k=1}^{t-1} r_k\nu_k^2  + \nu_t(\sum_{i=1}^n d_i- r_t \nu_t ) \big) \cr
&=& - (n-1) \nu_t^2 - \sum_{i=1}^n d_i^2
+  \sum_{k=1}^{t-1} r_k \nu_k^2 + \nu_t^2(n-1-r_t)
+2 r_t \nu_t^2   \cr
&=& \sum_{k=1}^t r_k \nu_k^2 - \sum_{i=1}^n d_i^2
\end{eqnarray*}
This gives the result.
\qed

\medskip
The rationality of the Hilbert-Kunz multiplicity
does not require the conditions normal and standard-graded,
as the following corollary shows.
Also the condition that the ground field is algebraically closed is not essential,
we only need that the ring is geometrically irreducible.

\begin{corollary}
\label{rationalcorollary} Let $R$ denote an $\NN$-graded
two-dimensional domain of finite type over an algebraically closed
field $K$ of positive characteristic. Let $I$ denote a homogeneous
$R_+$-primary ideal. Then the Hilbert-Kunz multiplicity $e_{HK}(I)$
is rational.
\end{corollary}

\proof By adjoining suitable roots for the algebra generators of $R$
we get a standard-graded $K$-domain $R \subseteq S$ finite over $R$
(see the proof of Theorem 4.2 in \cite{brennertightplus}). Due to
\cite[Theorem 2.7]{watanabeyoshidahilbertkunzinequality} we have for
finite extensions the relationship $e_{HK}(I) = e_{HK} (IS)/ \rk
(S)$, so we may assume that $R$ is a standard-graded domain. Its
normalization $R \subseteq \tilde{R}$ is a graded domain (see
\cite[\S 62, Aufgabe 27]{schejastorch2}), which might not be
standard-graded. However the open subsets $D_+(x)$ for elements $x
\in R_1 \subseteq \tilde{R}_1$ do cover $Y= \Proj \tilde{R}$, and,
as remarked in footnote 1, Theorem \ref{hilbertkunzformula} also
holds under this assumption. \qed

\section{Remarks and Examples}

We gather together some corollaries of our main result and make several remarks.

\begin{remark}
How does the denominator of the Hilbert-Kunz multiplicity look like?
The formula in Theorem \ref{hilbertkunzformula} shows that possible factors are
$2$, $\deg (Y)$,
the numbers $r < n$, where $n$ is the number of ideal generators,
and some powers of the characteristic $p$.
\end{remark}

\begin{remark}
Theorem \ref{strongsemistabletheorem} is of course a special case of Theorem \ref{hilbertkunzformula}.
If the syzygy bundle is strongly semistable,
then $t=1$, $r_t=n-1$, $\nu_t= \frac{d_1 \plusdots d_n}{n-1}$.
\end{remark}

The following corollary treats the next easiest case, namely $t=2$,
so that the syzygy bundle is an extension of two strongly semistable
bundles.

\begin{corollary}
\label{t=2}
Suppose the situation \ref{situation} and suppose that $t=2$.
Then the Hilbert-Kunz multiplicity $e_{HK}(I)$ is given by
$$\frac{ \deg (Y)}{2}
\big( r_2\nu_2^2 +\frac{(\sum_{i=1}^n d_i -r_2 \nu_2)^2}{n-1-r_2}  - \sum_{i=1}^n d_i^2
 \big) \, .$$
\end{corollary}
\proof
This follows directly from Theorem \ref{hilbertkunzformula} using
$r_1=n-1-r_2$ and $r_1\nu_1= \sum_{i=1}^n d_i -r_2\nu_2$.
\qed

\begin{corollary}
\label{n=3}
Suppose the situation \ref{situation} and suppose that $n=3$.
Then the syzygy bundle has rank two and the following
two cases may occur.

\numiii

\begin{enumerate}

\item
The syzygy bundle is strongly semistable.
Then the Hilbert-Kunz multiplicity is
$$e_{HK}(I)= \frac{\deg (Y)}{2} \big( \frac{(d_1+d_2+d_3)^2}{2}  -  d_1^2-d_2^2-d_3^2 \big) \, .$$

\item
The syzygy bundle is not strongly semistable.Then $\nu_2 > \nu_1$
and the Hilbert-Kunz multiplicity $e_{HK}(I)$ is given by
$$ \deg (Y) (  \nu_2^2 -  \nu_2 \sum_{i=1}^3 d_i + \sum_{i < j} d_i d_j) \, .$$
\end{enumerate}
\end{corollary}
\proof
The first statement follows from Theorem \ref{strongsemistabletheorem}
(or Theorem \ref{hilbertkunzformula}).
For the second statement
suppose that the $e$-th pull-back is not semistable.
Then there exist
invertible sheaves $ 0 \ra \shL \ra \Syz(f_1^q,f_2^q,f_3^q) \ra \shM \ra 0$,
$q =p^{e}$, with $\deg (\shL) > \deg (\shM)$.
Such a filtration is strong.
Therefore $\nu_1 = - \frac{\deg(\shL)}{q \deg (Y)}$ and $\nu_2 = - \frac{\deg(\shM)}{q \deg (Y)}$.
We insert $r_1=r_2=1$ and $\nu_1+\nu_2= d_1+d_2+d_3$
in the formula of Corollary \ref{t=2} and get up to the factor $\deg (Y)/2$ the expression
\begin{eqnarray*}
\nu_2^2 +(\sum_{i=1}^3 d_i - \nu_2)^2 - \sum_{i=1}^3 d_i^2
&=& 2 \nu_2^2 - 2 \nu_2 \sum_{i=1}^3 d_i + (\sum_{i=1}^3 d_i)^2 - \sum_{i=1}^3 d_i^2              \cr
&=& 2 \nu_2^2 - 2 \nu_2 \sum_{i=1}^3 d_i +2 \sum_{i < j} d_i d_j
\end{eqnarray*}
Multiplying with $\deg (Y)/2$ gives the result.
\qed

\begin{remark}
If in the previous corollary the degrees are equal, then the formula
in the second case reduces to $\frac{\deg(Y)}{2} (  \nu_2^2 - 3
\nu_2 d +3d^2)$. Corollary \ref{n=3} is independent of Langers
theorem, since in rank two it is clear that either $\shS$ is
strongly semistable or some Frobenius pull-back of it has an
invertible subsheaf as in the proof of \ref{n=3}. The same is true
for the following Corollary, which treats the Hilbert-Kunz
multiplicity of the cone over a plane curve.
\end{remark}

\begin{corollary}
\label{planecurve}
Let $C=V_+(H) \subset \PP^2$ denote a smooth plane projective curve of degree $h$, $R=K[X,Y,Z]/(H)$.
Then the following hold.

\numiii

\begin{enumerate}

\item
There exists a rational number $ \frac{3}{2} \leq  \nu_2  \leq 2$
such that
$e_{HK}(R) =h (\nu_2^2 -3 \nu_2 +3 )$.

\item
$ \nu_2 = 3/2$ holds if and only if
the restriction of the tangent bundle $\shT_{\PP^2}$
to the curve $C$ is strongly semistable.
In this case we have $e_{HK}(R) = \frac{3}{4}h$.

\item
We have estimates
$ \frac{3}{4}h \leq e_{HK}(R) \leq h $.

\end{enumerate}
\end{corollary}
\proof
(i),(ii).
Note that the tangent bundle is dual to the syzygy bundle
$\Syz(X,Y,Z)$ of the variables. If the restriction of this bundle is strongly
semistable, then (ii) holds due to Corollary \ref{n=3} (i).
Then also (i) is true for $\nu_2 = 3/2$.

The maximal slope of $\Syz(X,Y,Z)|C \subset \bigoplus_3 \O(-1)$
cannot exceed $-h$, hence the minimal slope is at least $-2h$.
Therefore $\nu_2 \leq 2$. The restriction is not strongly semistable
if and only if $\nu_2 > 3/2$ holds. In this case the formula in
Corollary \ref{n=3} (ii) gives $(i)$.

(iii). The quadratic polynomial takes its minimum at $\nu_2= 3/2$.
The value at $\nu_2=2$ is $1$. \qed

\begin{remark}
We discuss the estimates in Corollary \ref{planecurve} (iii)
and relate it to some results in the literature.
The bound  $3h/4 \leq e_{HK}(R)$ for cones over plane curves
was obtained in \cite{buchweitzchenhilbertkunz}.

In general there exist estimates
$\frac{e(I)}{d!} \leq e_{HK } (I) \leq e(I)$ (\cite[Lemma 6.1]{hunekeapplication}),
where $e(I)$ denotes the Hilbert-Samuel multiplicity of $I$ in a $d$-dimensional local ring.
Recall that the Hilbert-Samuel multiplicity
of the cone over a projective variety is the degree of the variety.
Hence these estimates yield in our situation $ h/2 \leq e_{HK}(R) \leq h$.

Watanabe and Yoshida have shown for a $2$-dimensional local
Cohen-Macaulay ring of positive characteristic that $e_{HK}(R)  \geq
\frac{e+1}{2}$ holds. We have $\frac{3}{4}h \geq \frac{h+1}{2}$, and
equality holds exactly for $h=1,2$.
\end{remark}

\begin{remark}
\label{hilbertkunzchar0} We take the expression in the formula in
Theorem \ref{hilbertkunzformula} as a definition for the
Hilbert-Kunz multiplicity in characteristic $0$. For the slopes we
just have to take the slopes in the Harder-Narasimhan filtration of
the syzygy bundle. This Hilbert-Kunz multiplicity is in fact
independent of the ideal generators, as the following proposition
shows.
\end{remark}

\begin{proposition}
Suppose that $R$ is a two-dimensional standard-graded normal domain over an algebraically closed
field $K$ of characteristic $0$. Let $I=(f_1 \komdots f_n)$ denote a homogeneous $R_+$-primary ideal.
Then the Hilbert-Kunz multiplicity $e_{HK}(f_1 \komdots f_n)$ is independent
of the ideal generators.
\end{proposition}
\proof
We write temporarily $e_{HK} (f_1 \komdots f_n)$ instead of $e_{HK}(I)$.
It is enough to show that
$$e_{HK}(f_1 \komdots f_n) = e_{HK}(f_1 \komdots f_n, f) \, ,$$
where $f$ is a homogeneous element $f \in (f_1 \komdots f_n)$.
Let
$\shS_1 \subsetdots \shS_t=\Syz(f_1 \komdots f_n)$
denote the Harder-Narasimhan filtration of $\Syz(f_1 \komdots f_n)$,
let $\mu_k= \mu( \shS_k/ \shS_{k-1})$, $r_k =\rk (\shS_k/ \shS_{k-1})$ and $ \nu_k = -\mu_k/\deg(Y)$.
These numbers determine
$e_{HK}(f_1 \komdots f_n)= \frac{\deg(Y)}{2} (\sum_{k=1}^t r_k\nu_k^2- \sum_{i=1}^n d_i^2)$.
Suppose that the degree of $f$ is $e$ and that $f= \sum_{i=1}^n a_if_i$.
Then we have the relationship
$$ \Syz (f_1 \komdots f_n,f) \cong  \Syz(f_1 \komdots f_n) \oplus \O(-e) \, .$$
The mappings from right to left
are given by $(s_1 \komdots s_n) \mapsto (s_1 \komdots s_n, 0)$
and $1 \mapsto (-a_1 \komdots -a_n ,1)$.
Let $i$ be such that $ \nu_i \leq e < \nu_{i+1}$ or equivalently that
$ \mu_i \geq \mu( \O(-e)) = -e /\deg(Y) > \mu_{i+1} $
(suppose in the following that $\nu_i < e$ holds, the case $=$ is similar).
Then the Harder-Narasimhan filtration of $\Syz(f_1 \komdots f_n, f)$ is
$$ \shS_1 \subsetdots \shS_i \subset \shS_i \oplus \O(-e) \subset \shS_{i+1} \oplus \O(-e)
\subsetdots \shS_t \oplus \O(-e) \, .$$
The (semistable) quotient sheaves are then
$$ \shS_1,\, \shS_2/\shS_1, \komdots \shS_i/\shS_{i-1}, \, \O(-e),\, \shS_{i+1}/\shS_i \komdots
\shS_{t}/\shS_{t-1} \, ,$$
and the new ranks are $r_1 \komdots r_i,1, r_{i+1} \komdots r_t$
and the new degree thresholds are
$\nu_1 \komdots \nu_i,e, \nu_{i+1} \komdots \nu_t$.
Hence $e_{HK} (f_1 \komdots f_n, f) = e_{HK}(f_1 \komdots f_n)$.
\qed

%\begin{remark}
%It is  both interesting and difficult to study in a relative setting the behavior
%of the (strong) Harder-Narasimhan filtration and of the Hilbert-Kunz multiplicity
%in dependence of the prime number.
%\end{remark}

\begin{remark}
It is of course tempting to conjecture a version of Theorem \ref{hilbertkunzformula}
in higher dimensions by replacing the degree of $\O_Y(1)$ by the top self
intersection number $(\O_Y(1))^{\dim Y}$.
However it is not clear whether the slopes carry enough information to
control the intermediate cohomologies $H^{i}(Y, \Syz)$, $0< i < \dim Y$.
\end{remark}

\begin{remark}
We briefly explain the relation of our main result to tight closure.
In \cite{brennertightplus}
we proved using the strong Harder-Narasimhan filtration
that the tight closure and the plus closure coincide
for a homogeneous ideal
in a two-dimensional graded domain
of finite type over the algebraic closure of a finite field.
In fact we showed that the containment $f \in (f_1 \komdots f_n)^*$
is a property of the corresponding cohomology class
$c= \delta(f)  \in H^1(Y, \Syz(f_1 \komdots f_n)(m))$ ($m= \deg (f)$)
in the strong Harder-Narasimhan filtration.

Since it is known (see \cite[Theorem 5.4]{hunekeapplication}) that
in an analytically unramified and formally equidimensional local
ring $R$ the equation $e_{HK}(I) = e_{HK} (J) $ holds if and only if
$I^* =J^*$ holds true for two ideals $I \subseteq J$ which are
primary to the maximal ideal, it is not surprising that the behavior
of the Hilbert-Kunz function and the Hilbert-Kunz multiplicity is
encoded in the strong Harder-Narasimhan filtration of a syzygy
bundle of ideal generators.

In \cite{brennerhilbertkunzcriterion} we show that in characteristic
zero the Hilbert-Kunz multiplicity as defined in Remark
\ref{hilbertkunzchar0} has the same relationship to solid closure as
the Hilbert-Kunz multiplicity in positive characteristic has to
tight closure.
\end{remark}

\begin{example}
Consider the monomial ideal
$I=(X^3,XY^2,ZY^2)$ in $R=K[X,Y,Z]/(H)$, where $H$ is a homogeneous polynomial of degree $h$
such that $R$ is normal and such that $Z$ and $X$ are parameters.
Denote by $C= \Proj R$ the corresponding smooth projective curve.
The tripel $(0,Z,-X)$ defines a global syzygy of degree $4$ (without common zero)
and so we get the short exact sequence
$$0 \lra \O(-4) \lra \Syz(X^3,XY^2,ZY^2)  \lra \O(-5) \lra 0 \, ,$$
and the inclusion $\O(-4) \subset \Syz(X^3,XY^2,ZY^2) $ is the Harder-Narasimhan filtration.
This filtration is of course strong, and the slope numbers
are $\mu_1 =-4h$ and $\mu_2= -5h$ (hence $\nu_1=4$ and $\nu_2=5$)
and the ranks are $r_1=r_2=1$.
The formula in Corollary \ref{t=2} (ii) yields that the Hilbert-Kunz
multiplicity is
$h(25-5 \cdot 9 +3 \cdot 9)=7h$.
\end{example}

\bibliographystyle{plain}

\bibliography{bibliothek}

\end{document}